# A stability preserved time-integration method for nonlinear advection-diffusion models

Huseyin Tunc[a], Murat Sari[b]

Department of Mathematics, Faculty of Arts and Science, Yildiz Technical University, Istanbul 34220, Turkey.

[a] tnchsyn@gmail.com, ORCID: 0000-0001-6450-5380

[b] sarim@yildiz.edu.tr, ORCID: 0000-0003-0508-2917

**Abstract**

A new implicit-explicit local differential transform method (IELDTM) is derived here for time integration of the nonlinear advection-diffusion processes represented by (2+1)-dimensional Burgers equation. The IELDTM is adaptively constructed as stability preserved and high order time integrator for spatially discretized Burgers equation. For spatial discretization of the model equation, the Chebyshev spectral collocation method (ChCM) is utilized. A robust stability analysis and global error analysis of the IELDTM are presented with respect to the direction parameter $\theta$. With the help of the global error analysis, adaptivity equations are derived to minimize the computational costs of the algorithms. The produced method is shown to eliminate the accuracy disadvantage of the classical $\theta$-method and the stability disadvantages of the existing DTM-based methods. Two examples of the Burgers equation in one and two dimensions have been solved via the ChCM-IELDTM hybridization, and the produced results are compared with the literature. The present time integrator has been proven to produce more accurate numerical results than the MATLAB solvers, *ode45* and *ode15s*.

*Keywords*: Local transform method; Taylor series; Advection-diffusion processes; Stability analysis

MSC2010: 65L05, 65L06, 34K28

## 1. Introduction

Time integrations of the reduced space-time partial differential equations are of great importance. Various versions of both explicit and implicit time integrators are always under development to create both accurate and stable numerical solvers. For higher dimensional parabolic PDEs, the computational efficiency of the time integrator becomes more critical due to the large ODE systems obtained after spatial discretization. In such cases, performance of the effective explicit ODE solvers such as Taylor series methods or explicit Runge-Kutta methods is weakened due to instability [1-3]. Implicit time integration schemes are the best options in such situations, but require more computational effort, especially for nonlinear problems [1-3].

Time integration methods are mainly divided into two subgroups, such as multi-point [4-9] and multi-derivative methods [10-18]. The multi-point methods are based on interpolation, such as finite element and spectral methods. Even if this group of methods is widely used for the spatial parts of PDEs, it is possible to solve also IVPs effectively with these methods. The multipoint methods are stability preserved and even capable of solving stiff IVPs. Despite all these flexibilities, this group of methods leads to a large degrees of freedom and a large system of linear/nonlinear algebraic equations for solving ODEs due to their non-iterative nature. This is an essential drawback due to both computational cost and loss of accuracy over large time scales. The multi-derivative methods such as Runge-Kutta methods [10-12], Rosenbrock methods [13-15], Taylor series methods [16-18] have been used extensively to observe the behavior of stiff and non-stiff IVPs. As for this group of methods, they are iterative in nature and have both explicit and implicit forms of their own. The Runge-Kutta and Rosenbrock methods are derived from the Taylor series expansions by introducing interior ghost nodes to get higher-order approximations. The Taylor series methods directly use the higher order derivatives of given implicit/explicit functions by the procedure known as automatic differentiation or differential transformation [16-23]. The idea of differential transformation leads us to a useful procedure for recursively evaluating higher order Taylor coefficients [19-23]. With the work of Tunc and Sari [24], a significant progress was made to solve stiff IVPs with the idea of differential transformation. The derived method is called the implicit/explicit local differential transformation method (IELDTM), which is based on the combination of the idea of differential transformation and the Taylor series method.

When dealing with the time integration of the advection-diffusion processes, numerous numerical methods have been proposed in the literature [25–36]. Explicit Runge-Kutta schemes are more common in applications [25-29] due to their efficient computational structures. It is well known that this group of time-integrators has lack of stability preserving properties. The second-order and unconditionally stable Crank-Nicolson method (CNM) is another widely used time integrator in advection-diffusion processes, and the CNM is a special case of the well-known θ-method [30-34]. Even if the stability properties of the CNM was found to be excellent, the method is second order and needs small time-increments to get higher accuracy. Here the CNM has proven to be a special case of the current IELDTM and therefore the restriction of the Crank-Nicolson approach has been overcome by the IELDTM by increasing its order. Another time-integration technique for solving reduced ODEs is the backward differentiation method (BDM), in which stability is preserved and has higher-order formulae [35–36]. The BDM has $A$ −stable formulae up to the second-order accuracy, while the present IELDTM has $A$ −stable formulae up to the fourth-order accuracy [24].

Here we have derived an implicit-explicit local differential transform method for time integration of the (2+1)-dimensional Burgers equation. The Chebyshev spectral collocation method (ChCM) has been used for the spatial discretization of considered PDEs, due to the important interpolation characteristics of the method [27]. The flexibility, stability characteristics, adaptivity and computational efficiency of the IELDTM [24] led us to improve our work in solving nonlinear advection-diffusion equations with the ChCM-IELDTM hybridization. It has been shown that the current hybrid approach destroys the disadvantageous of the existing numerical methods such as $\theta$ − methods and DTM-based methods. A priori error analysis has been done and convergence rates are determined according to the direction parameter $\theta$. Numerical properties of the current hybrid approach have been compared with the literature and MATLAB solvers, *ode15s* and *ode45*. It has been shown that the present approach offers numerically better performance than the previously produced methods, and the IELDTM is found to be a great option for time integrations of the PDEs.

## 2. ChCM-IELDTM Hybridization

In this section, we introduce the implementation procedure of the IELDTM-ChCM hybridization to 2D Burgers equation represents the nonlinear advection-diffusion mechanism. Consider the following 2D Burgers equation

$$u_t + uu_x + uu_y = \varepsilon(u_{xx} + u_{yy}), \quad (x,y) \in D, t > 0 \tag{1}$$

with the boundary conditions

$$u(a,y,t) = f_1(y,t), \quad u(b,y,t) = f_2(y,t) \tag{2}$$

$$u(x,c,t) = h_1(x,t), \quad u(x,d,t) = h_2(x,t) \tag{3}$$

and initial condition

$$u(x,y,0) = g(x,y) \tag{4}$$

where $D = \{(x,y)|\ a \leq x \leq b, c \leq y \leq d\}$, $\varepsilon$ is kinematic viscosity constant for $\varepsilon > 0$ and $f_1$, $f_2$, $h_1$, $h_2$, and $g$ are known smooth functions. The subscripts $x, y$ and $t$ represent differentiations with respect to spaces $x$, $y$ and time $t$, respectively.

Let us approximate the solution $u(x, y, t)$ with $(N, M)$ order Chebyshev polynomials as follows [27]

$$u(x, y, t) = \sum_{n=0}^{N} \sum_{m=0}^{M} \alpha_{nm} \beta_{nm}(t) \bar{T}_n(x) \bar{T}_m(y) \qquad (5)$$

where $\beta_{nm}(t)$ are time dependent parts of the representation, $\alpha_{nm} = 1$ for interior points and $\alpha_{nm} = 1/2$ for any boundary point. In equation (5), $\bar{T}_n(x)$ and $\bar{T}_m(y)$ are defined in terms of the first kind Chebyshev polynomials $T_n(x)$ and $T_m(y)$ as

$$\bar{T}_n(x) = T_n\left(\frac{2x-(a+b)}{b-a}\right) = \cos\left(n \arccos\left(\frac{2x-(a+b)}{b-a}\right)\right), \qquad (6)$$

$$\bar{T}_m(y) = T_m\left(\frac{2y-(c+d)}{d-c}\right) = \cos\left(m \arccos\left(\frac{2y-(c+d)}{d-c}\right)\right). \qquad (7)$$

We define the following restricted collocation points on $[-1, 1]$ as

$$x_n = \frac{1}{2}\left((a+b) - (b-a)\cos\left(\frac{\pi n}{N}\right)\right), \quad n = 0,1,2,\ldots,N, \qquad (8)$$

$$y_m = \frac{1}{2}\left((c+d) - (d-c)\cos\left(\frac{\pi m}{M}\right)\right), \quad m = 0,1,2,\ldots,M.$$

The discrete orthogonality relation is vital for the ChCM defined as

$$\sum_{n=0}^{N} \alpha_n \bar{T}_i(x_n) \bar{T}_j(x_n) = \gamma_i \delta_{ij} \qquad (9)$$

$$\sum_{m=0}^{M} \alpha_m \bar{T}_i(y_m) \bar{T}_j(y_m) = \omega_i \delta_{ij} \qquad (10)$$

with

$$\gamma_i = \begin{cases} \frac{N}{2}, & i \neq 0, N \\ N, & i = 0, N, \end{cases} \qquad (11)$$

$$\omega_i = \begin{cases} \frac{M}{2}, & i \neq 0, M \\ M, & i = 0, M. \end{cases} \qquad (12)$$

Thus, we reach the required derivatives at predetermined collocation points in equation (1) as follows:

$$u_x(x_i, y_j, t) = \sum_{n=0}^{N} [A_x]_{in} u_{nj}(t), \qquad (13)$$

$$u_{xx}(x_i, y_j, t) = \sum_{n=0}^{N} [B_x]_{in} u_{nj}(t), \qquad (14)$$

$$u_y(x_i, y_j, t) = \sum_{m=0}^{M} [A_y]_{jm} u_{im}(t), \qquad (15)$$

$$u_{yy}(x_i, y_j, t) = \sum_{m=0}^{M} [B_y]_{jm} u_{im}(t), \qquad (16)$$

where $i = 0,1,2,\ldots,N$, $j = 0,1,2,\ldots,M$, $u_{nm}(t) = u(x_n, y_m, t)$, $B_x = A_x^2$, $B_y = A_y^2$, $A_x$ and $A_y$ are $(N \times N)$ and $(M \times M)$ matrices are defined by

$$[A_x]_{in} = \begin{cases} \mu_1 \sum_{j=0}^{N}(-1)^j j^2 \cos\left(j\left(\pi - \frac{(N+n-1)\pi}{N}\right)\right), & i = 0 \\ \mu_1 \sum_{j=0}^{N} \frac{j\sin\left(j\left(\pi - \frac{(i-1)\pi}{N}\right)\right)\cos\left(j\left(\pi - \frac{(n-1)\pi}{N}\right)\right)}{\sqrt{1-\cos^2\left(\pi - \frac{(i-1)\pi}{N}\right)}}, & i = 1, 2, \ldots, N-1 \\ \mu_1 \sum_{j=0}^{N}(-1)^{j+1} j^2 \cos\left(j\left(\pi - \frac{(n-1)\pi}{N}\right)\right), & i = N, \end{cases} \quad (17)$$

$$[A_y]_{im} = \begin{cases} \mu_2 \sum_{j=0}^{M}(-1)^j j^2 \cos\left(j\left(\pi - \frac{(M+m-1)\pi}{M}\right)\right), & i = 0 \\ \mu_2 \sum_{j=0}^{M} \frac{j\sin\left(j\left(\pi - \frac{(i-1)\pi}{M}\right)\right)\cos\left(j\left(\pi - \frac{(m-1)\pi}{M}\right)\right)}{\sqrt{1-\cos^2\left(\pi - \frac{(i-1)\pi}{N}\right)}}, & i = 1, 2, \ldots, M-1 \\ \mu_2 \sum_{j=0}^{M}(-1)^{j+1} j^2 \cos\left(j\left(\pi - \frac{(m-1)\pi}{M}\right)\right), & i = M \end{cases} \quad (18)$$

where $\mu_1 = \frac{2}{b-a}$, $\mu_2 = \frac{2}{d-c}$, $[A_x]_{in}$ and $[A_y]_{im}$ are called the Chebyshev differentiation matrices. Writing equations (16)-(19) into equation (1) and imposing the boundary conditions leads to the following nonlinear ODE system,

$$\frac{d\boldsymbol{\beta}}{dt} = \varepsilon(\bar{B}_x \boldsymbol{\beta} + \boldsymbol{\beta}\bar{B}_y^T) - \langle \boldsymbol{\beta}, (\bar{A}_x \boldsymbol{\beta} + \boldsymbol{\beta}\bar{A}_y^T)\rangle + \boldsymbol{F}(t) \quad (19)$$

where $\langle .,. \rangle$ is the standard elementwise product, $\boldsymbol{\beta}(t) = [\boldsymbol{\beta}(t)]_{nm}$ is $((N-1) \times (M-1))$ is matrix of time variable, $\boldsymbol{F}(t) = [\boldsymbol{F}(t)]_{nm}$ is $((N-1) \times (M-1))$ matrix of time variable related to boundary conditions (2)-(3). After applying the boundary conditions to the matrices $A_x, A_y, B_x$ and $B_y$, respectively, $\bar{A}_x, \bar{A}_y, \bar{B}_x,$ and $\bar{B}_y$ are obtained. Note that the initial condition of the nonlinear ODE system is defined by

$$\boldsymbol{\beta}(0) = [\boldsymbol{\beta}(0)]_{nm} = g(x_n, y_m). \quad (20)$$

Now, we construct the IELDTM to solve IVP stated in (18)-(19). Detailed description of the IELDTM for the solution of IVPs can be found in the literature [24]. Let us divide the time interval $[0, t_f]$ into at most $P$ time elements with $\Delta t_i = t_{i+1} - t_i$ and the partition of the interval as $\omega = \{0 = t_0 < t_1 < \cdots < t_{P^*} = t_f\}$. Let us then consider the convergent Taylor series representation of the function $\boldsymbol{\beta}(t)$ of order $K$ about $t = t_i$ as

$$\boldsymbol{\beta}_i(t) = \sum_{k=0}^{K} \bar{\boldsymbol{\beta}}_i(k)(t-t_i)^k + O((t-t_i)^{K+1}), \quad t_i - \rho^i \leq t \leq t_i + \rho^i \quad (21)$$

where $i = 0, 1, \ldots, P^*$, $\bar{\boldsymbol{\beta}}_i(k)$ is the local differential transformation of $\boldsymbol{\beta}_i(t)$ and $\rho^i$ is the radius of convergence of the representation. Taking differential transform of equation (19) leads to,

$$\bar{\boldsymbol{\beta}}_i(k+1) = \frac{1}{k+1}\left[\varepsilon(\bar{B}_x \bar{\boldsymbol{\beta}}_i(k) + \bar{\boldsymbol{\beta}}_i(k)\bar{B}_y^T) - \sum_{s=0}^{k} \bar{\boldsymbol{\beta}}_i(s)(\bar{A}_x \bar{\boldsymbol{\beta}}_i(k-s) + \bar{\boldsymbol{\beta}}_i(k-s)\bar{A}_y^T) + \bar{\boldsymbol{F}}_i(k)\right] \quad (22)$$

where $\bar{\boldsymbol{F}}_i(k)$ is the local differential transformation of $\boldsymbol{F}(t)$. By defining the parameter $\theta \in [0,1]$ as a direction parameter, continuity requirement of the two consecutive solutions yields the following equation system,

$$\sum_{k=0}^{K} \overline{\beta}_{i+1}(k)(-\theta \Delta t_i)^k = \sum_{k=0}^{K} \overline{\beta}_i(k)((1-\theta)\Delta t_i)^k + O((\Delta t_i)^{K+1}, \theta) \tag{23}$$

where $O((\Delta t_i)^{K+1}, \theta)$ represents the dependency of the local truncation error to time increment, transformation order and direction parameter. By defining $\overline{\beta}_0(0) = \beta(0)$, the rest of $\overline{\beta}_i(k)$ can be calculated from the recursive relation (22). Since $\overline{\beta}_i(0)$ is known from the previous step, equation (23) is an implicit/explicit algebraic equation system of $\overline{\beta}_{i+1}(0)$. This numerical algorithm finds all local solutions $\beta_i(t)$ by solving equation (23) at each time step. Finally, by obtaining all local solutions $\beta_i(t)$ we reach

$$u^i(x, y, t) = \sum_{n=0}^{N} \sum_{m=0}^{M} \alpha_{nm} (\beta_{nm})_i(t) \overline{T}_n(x) \overline{T}_m(y) \tag{24}$$

where $i = 0, 1, \ldots, P^*$.

## 3. Error Analysis

In this section, priori error estimation of the IELDTM for solving the nonlinear ODE system (19) will be presented. Let us first define the map $\phi: R^{(N-1) \times (M-1)} \to R^{(N-1)(M-1) \times 1}$ as

$$\phi(\beta_i(k)) = \left[(\beta_i(k))_{11}, (\beta_i(k))_{12}, \ldots, (\beta_i(k))_{21}, \ldots, (\beta_i(k))_{(N-1)(M-1)}\right]^T = \varphi_i(k)$$

The map $\phi$ reshapes the matrix $\beta_i(k)$ as a related vector $\varphi_i(k)$ for further analysis. The exact form of the continuity equation (23) can be written as

$$\overline{\overline{\varphi}}_{i+1}(0) = \overline{\overline{\varphi}}_i(0) + \sum_{k=0}^{K} \overline{\overline{\varphi}}_i(k)((1-\theta)\Delta t)^k - \sum_{k=0}^{K} \overline{\overline{\varphi}}_{i+1}(k)(-\theta \Delta t)^k + \Delta t \rho_i \tag{25}$$

where $\overline{\overline{\varphi}}$ is the exact form of $\overline{\varphi}$ and then the local truncation error $\rho_i$ takes the following form

$$\rho_i = \left[\overline{\overline{\varphi}}_{(i+1)_*}(K+1)(1-\theta)^{K+1} - \overline{\overline{\varphi}}_{i_*}(K+1)(-\theta)^{K+1}\right] dx^K. \tag{26}$$

Here $\overline{\overline{\varphi}}_{i_*}$ defines the local transform at the local point $t_{i_*} \in [t_{i-1}, t_{i+1}]$ for all $i = 0,1,2,\ldots,P-1$. With the use of recursive relation (22), equation (25) can be rewritten as

$$\overline{\overline{\varphi}}_{i+1}(0) = \overline{\overline{\varphi}}_i(0) + (1-\theta)\Delta t H_1(\overline{\overline{\varphi}}_i(0), \theta, \Delta t) + \theta \Delta t H_2(\overline{\overline{\varphi}}_{i+1}(0), \theta, \Delta t) + \Delta t \rho_i \tag{27}$$

where $H_1: R^{(N-1)(M-1) \times 1} \times [0,1] \times R^+ \to R^{(N-1)(M-1) \times 1}$ and $H_2: R^{(N-1)(M-1) \times 1} \times [0,1] \times R^+ \to R^{(N-1)(M-1) \times 1}$ are defined by

$$H_1(\overline{\overline{\varphi}}_i(0), \theta, \Delta t) = \sum_{k=1}^{K} \overline{\overline{\varphi}}_i(k)((1-\theta)\Delta t)^{k-1}, \tag{28}$$

$$H_2(\overline{\overline{\varphi}}_{i+1}(0), \theta, \Delta t) = \sum_{k=1}^{K} \overline{\overline{\varphi}}_{i+1}(k)(-\theta \Delta t)^{k-1}. \tag{29}$$

The IELDTM form of the continuity equation (23) takes the form

$$\overline{\varphi}_{i+1}(0) = \overline{\varphi}_i(0) + (1-\theta)\Delta t H_1(\overline{\varphi}_i(0), \theta, \Delta t) + \theta \Delta t H_2(\overline{\varphi}_{i+1}(0), \theta, \Delta t). \tag{30}$$

Let us define the following Jacobian matrices

$$J_1(\varphi, \theta, \Delta t) = \frac{\partial (H_1)_i}{\partial \varphi_j}, \tag{31}$$

$$J_2(\varphi, \theta, \Delta t) = \frac{\partial (H_2)_i}{\partial \varphi_j}. \tag{32}$$

Assume that the conditions

$$\mu(J_1(\varphi,\theta,\Delta t)) \leq 0 \text{ and } \mu(J_2(\varphi,\theta,\Delta t)) \leq 0 \tag{33}$$

hold for all $(\varphi,\theta,\Delta t) \in R^{(N-1)(M-1)\times 1} \times [0,1] \times R^+$ and $\mu(A) = \lim_{\tau\to 0}\frac{\|I+\tau A\|-1}{\tau}$ is the logarithmic norm [39]. Let $Z_i^1 \in R^{(N-1)(M-1)\times(N-1)(M-1)}$ and $Z_i^2 \in R^{(N-1)(M-1)\times(N-1)(M-1)}$ be defined for $i \geq 0$ with

$$Z_i^1 = \int_0^1 \Delta t J_1(\sigma\bar{\bar{\varphi}}_i + (1-\sigma)\bar{\varphi}_i, \theta, \Delta t)d\sigma, \tag{34}$$

$$Z_{i+1}^2 = \int_0^1 \Delta t J_2(\sigma\bar{\bar{\varphi}}_{i+1} + (1-\sigma)\bar{\varphi}_{i+1}, \theta, \Delta t)d\sigma. \tag{35}$$

Applying the mean value theorem to equations (34)-(35) gives

$$\Delta t H_1(\bar{\bar{\varphi}}_i(0),\theta,\Delta t) - \Delta t H_1(\bar{\varphi}_i(0),\theta,\Delta t) = Z_i^1(\bar{\bar{\varphi}}_i(0) - \bar{\varphi}_i(0)), \tag{36}$$

$$\Delta t H_2(\bar{\bar{\varphi}}_{i+1}(0),\theta,\Delta t) - \Delta t H_2(\bar{\varphi}_{i+1}(0),\theta,\Delta t) = Z_{i+1}^2(\bar{\bar{\varphi}}_{i+1}(0) - \varphi_{i+1}(0)). \tag{37}$$

Defining global discretization error $\varepsilon_i = \varphi(t_i) - \bar{\varphi}_i(0) = \bar{\bar{\varphi}}_i(0) - \bar{\varphi}_i(0)$ and subtracting equation (30) from equation (27) give rise to

$$\varepsilon_{i+1} = \varepsilon_i + (1-\theta)Z_i^1\varepsilon_i + \theta Z_{i+1}^2\varepsilon_{i+1} + \Delta t\rho_i. \tag{38}$$

From (34)-(35) and Lemma 2.6 in literature [39], the following inequalities hold

$$\mu(Z_i^1) \leq 0 \text{ and } \mu(Z_i^2) \leq 0.$$

The analysis are considered in two main cases as follows:

**Case 1:** Assume that $\theta = 1$, i.e. a backward local differential transform method is obtained. Equation (38) gives

$$\varepsilon_{i+1} = (I - Z_{i+1}^2)^{-1}(\varepsilon_i + \Delta t\rho_i). \tag{39}$$

Since $\|(I - Z_{i+1}^2)^{-1}\| \leq 1$ according to equation (38) and Theorem 2.13 of the literature [39], we find

$$\|\varepsilon_i\| \leq \|\varepsilon_0\| + \sum_{j=0}^{i-1}\Delta t\|\rho_j\| \leq \|\varepsilon_0\| + t_i \max_{0\leq j\leq i-1}\|\rho_j\|. \tag{40}$$

Assumption of $\varepsilon_0 \to 0$ leads to the following global error estimate

$$\|\varphi(t_i) - \bar{\varphi}_i(0)\|_\infty = \|\beta(t_i) - \bar{\beta}_i(0)\|_\infty \leq t_i\omega\Delta t^K \tag{41}$$

where

$$\omega = \left|\frac{1}{(K+1)!} \max_{t\in[0,T]}\left(\left\|\frac{d\varphi^{K+1}}{dt^{K+1}}\right\|\right)\right|.$$

**Case 2:** Let us now consider the rest of the implicit cases for $\theta \geq 1/2$. The transformed errors are defined as

$$\bar{\varepsilon}_i = (I - \theta Z_i^2)\varepsilon_i \tag{42}$$

where $i = 0,1,\ldots,P$. Then the following recursive relation can be obtained from equation (37)

$$\bar{\varepsilon}_{i+1} = (I + (1-\theta)Z_i^1)(I - \theta Z_i^2)^{-1}\bar{\varepsilon}_i + \Delta t \rho_i = R(Z_i^1, Z_i^2)\bar{\varepsilon}_i + \Delta t \rho_i. \tag{43}$$

Assume also that the following stability condition holds

$$\|R(Z_i^1, Z_i^2)\| = \|(I + (1-\theta)Z_i^1)(I - \theta Z_i^2)^{-1}\| \leq 1. \tag{44}$$

Thus, equation (43) becomes

$$\|\bar{\varepsilon}_i\| \leq \|\bar{\varepsilon}_0\| + \sum_{j=0}^{i-1} \Delta t \|\rho_j\| \leq \|\bar{\varepsilon}_0\| + t_i \max_{0 \leq j \leq i-1} \|\rho_j\|. \tag{45}$$

Since $\|(I - \theta Z_i^2)^{-1}\| \leq 1$ for $\theta \geq 1/2$, the main global error estimate can be constructed from equation (45) as

$$\|\varepsilon_i\| \leq \|\varepsilon_0\| + t_i \max_{0 \leq j \leq i-1} \|\rho_j\|. \tag{46}$$

Assumption of $\varepsilon_0 \to 0$ leads to the following global error estimate

$$\|\varphi(t_i) - \bar{\varphi}_i(0)\|_\infty = \|\beta(t_i) - \bar{\beta}_i(0)\|_\infty \leq \begin{cases} t_i \omega \Delta t^{K+1}, & \text{if } \theta = 1/2 \text{ and } K \text{ is odd} \\ t_i \omega \Delta t^K, & \text{otherwise} \end{cases} \tag{47}$$

where

$$\omega = \begin{cases} \left|\frac{1}{(K+2)!} \max_{t \in [0,T]} \left(\left\|\frac{d\varphi^{K+2}}{dt^{K+2}}\right\|\right)\left(\frac{1}{2}\right)^{K+1}(K+1)\right|, & \text{if } \theta = 1/2 \text{ and } K \text{ is odd} \\ \left|\frac{1}{(K+1)!} \max_{t \in [0,T]} \left(\left\|\frac{d\varphi^{K+1}}{dt^{K+1}}\right\|\right)[(1-\theta)^K - (-\theta)^K]\right|, & \text{if } \theta \neq \frac{1}{2} \text{ and } K \text{ is arbitrary.} \end{cases} \tag{48}$$

Adaptive construction of the step-sizes is important to estimate the accuracy of the solution and to reduce the computational cost. With the help of the local truncation error defined in equation (47), adaptive step-sizes of the algorithm can be determined with the following inequalities

$$\Delta t_i < \left(\frac{tol}{\left|\left(\left(\frac{1}{2}\right)^{K+1}(K+1)\right)\|\bar{\beta}_i(K+2)\|_\infty\right|}\right)^{\frac{1}{K+1}}, \text{ if } \theta = 0.5 \text{ and } K \text{ is odd} \tag{49}$$

$$\Delta t_i < \left(\frac{tol}{|(1-\theta)^{K+1} - (-\theta)^{K+1}|\|\bar{\beta}_i(K+1)\|_\infty}\right)^{\frac{1}{K}}, \text{ otherwise} \tag{50}$$

where $K$ is the order of local polynomial approximation (21), $tol$ is the predetermined tolerance and $\bar{\beta}_i(K+1)$ or $\bar{\beta}_i(K+2)$ can be obtained from the recursive relation (22).

## 4. Numerical Experiments

In this section, we have numerically analyzed the present ChCM-IELDTM approach over the 1D-2D Burgers equations. The exponential convergence of the IELDTM has been proven for various values of the direction parameter $\theta$. The computational efficiency of the adaptive IELDTM has been shown in comparison with the literature and MATLAB solvers, *ode15s* and *ode45*.

**Problem 1 [31]**

Consider the one-dimensional form of the Burgers equation (1) with the initial condition

$$u(x,0) = \sin \pi x, \quad 0 < x < 1 \tag{51}$$

and the homogenous Dirichlet boundary conditions

$$u(0,t) = u(1,t) = 0, \quad t > 0. \tag{52}$$

The exact solution of the 1D form of the Burgers equation under the consideration of the conditions (51)-(52) is given by

$$u(x,t) = 2\pi\varepsilon \frac{\sum_{n=1}^{\infty} a_n \exp(-n^2\pi^2\varepsilon t) n \sin(n\pi x)}{a_0 + \sum_{n=1}^{\infty} a_n \exp(-n^2\pi^2\varepsilon t) \cos(n\pi x)} \tag{53}$$

with the Fourier coefficients

$$a_0 = \int_0^1 \exp\{-(2\pi\varepsilon)^{-1}[1 - \cos(\pi x)]\} dx, \tag{54}$$

$$a_n = 2\int_0^1 \exp\{-(2\pi\varepsilon)^{-1}[1 - \cos(\pi x)]\} \cos(n\pi x) dx. \tag{55}$$

The performance of the current ChCM-IELDTM is shown comparing the maximum errors depending on the varying values of the direction parameter θ (see Figure 1). As theoretically expected, the maximum errors are reduced exponentially by increasing the order of the method. It can be observed from Figure 1 that choosing $\theta = 0.5$ with odd number $K$ values leads to the numerical method of $K + 1$ order. Note that the exponential convergence is verified for both $\varepsilon = 1$ and $\varepsilon = 0.1$ selections. The advection-dominated cases of the present example exhibit challenging behavior as has been widely discussed in the literature [25,31,33]. The ChCM-IELDTM hybridization produces non-oscillatory solutions in such stiff cases as shown in Figure 2. The sharp behaviors are demonstrated by choosing the kinematic viscosity constants as $\varepsilon = 0.001$ and $\varepsilon = 0.0005$. Using the current central time integration method and high-degree spectral approach, challenging numerical solutions have been captured for both cases. The performance of time adaptive ChCM-IELDTM hybridization have been illustrated based on the error norms and the required step sizes as seen in Table 1. It is obvious from the table that the IELDTM time-integration is more preferable than the *ode45-ode15s* of MATLAB. The overall performance of the current approach is compared with finite element-based techniques such as the weak Galerkin finite element method (WFEM) [37] and the strong Galerkin finite element method (GFEM) [38] as shown in Table 2. In the temporal part of FEM based algorithms, the Crank-Nicolson method (a special case of IELDTM; = 0.5 and K = 1) has been used as a time integration technique. The current approach, ChCM-IELDTM, produces more reliable results than the FEM-based techniques, even if their degrees of freedom are lower than the literature [37–38]. The effect of the IELDTM is easily seen from the results in the table, i.e. the case of $K = 3$ and $\theta = 0.5$ (fourth order in time) reduces the maximum error from $5.57E - 07$ to $8.97E - 10$.

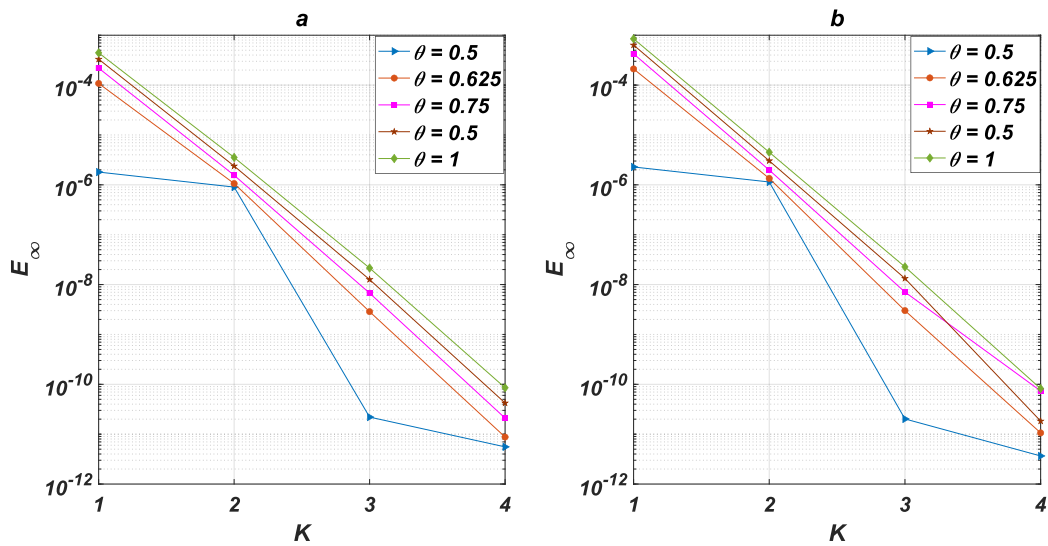

**Fig. 1** Order refinement results of the ChCM-IELDTM at $t_f = 0.5$ for various values of the direction parameter $\theta$ with the parameter values a) $\varepsilon = 1, N = 10$ and $dt = 0.0025$ b) $\varepsilon = 0.1, N = 40$ and $dt = 0.0025$.

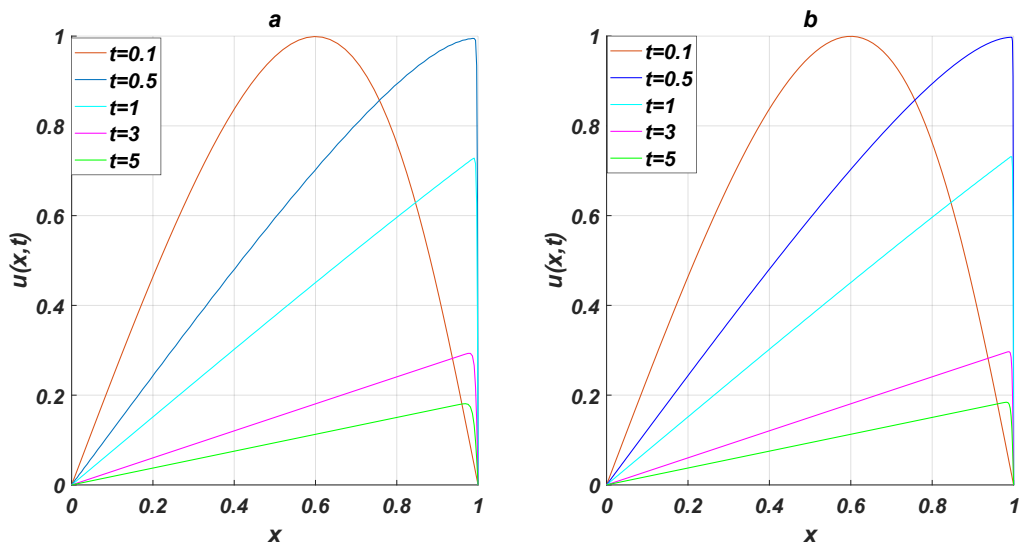

**Fig. 2** The ChCM-IELDTM solutions of Problem 1 at various times for the parameter values a) $\varepsilon = 0.001, N = 80, \theta = 0.5, K = 4$ and $dt = 0.001$ b) $\varepsilon = 0.0005, N = 120, \theta = 0.5, K = 4$ and $dt = 0.001$.

**Table 1** Number of time steps and the maximum pointwise errors of the time-adaptive ChCM-IELDTM approach for various values of the viscosity constant $\varepsilon$ and time approximation order $K$ with $\theta = 0.5$, $t_f = 0.5$ and $tol = E - 10$.

| $\varepsilon$ | Method | Step | $\|E\|_\infty$ |
|---|---|---|---|
| $\varepsilon = 1$ $N = 10$ | IELDTM-ChCSM $K = 3$ | 500 | 8.5E-09 |
| | IELDTM-ChCSM $K = 4$ | 354 | 8.5E-09 |
| | IELDTM-ChCSM $K = 5$ | 101 | 8.5E-09 |
| | *ode45*-ChCSM | 493 | 8.0E-09 |
| | *ode15s*-ChCSM | 243 | 8.6E-09 |
| $\varepsilon = 0.1$ $N = 20$ | IELDTM-ChCSM $K = 3$ | 266 | 2.0E-08 |
| | IELDTM-ChCSM $K = 4$ | 189 | 2.0E-08 |
| | IELDTM-ChCSM $K = 5$ | 88 | 2.0E-08 |
| | *ode45*-ChCSM | 596 | 2.0E-08 |
| | *ode15s*-ChCSM | 163 | 2.0E-08 |
| $\varepsilon = 0.01$ $N = 40$ | IELDTM-ChCSM $K = 3$ | 583 | 1.2E-05 |
| | IELDTM-ChCSM $K = 4$ | 413 | 1.2E-05 |
| | IELDTM-ChCSM $K = 5$ | 126 | 1.2E-05 |
| | *ode45*-ChCSM | 839 | 1.2E-05 |
| | *ode15s*-ChCSM | 253 | 1.2E-05 |

**Table 2** Comparison of the numerical results produced with $t_f = 0.1$, $\varepsilon = 0.1$ and $dt = 0.001$.

| $x$ | WGFEM [37] $N = 40$ | GFEM [38] $N = 40$ | ChCM-IELDTM ($\theta = 0.5 / K = 1$) $N = 20$ | ChCM-IELDTM ($\theta = 0.5 / K = 3$) $N = 20$ | Exact |
|---|---|---|---|---|---|
| $x = 0.1$ | 0.22346 | 0.2234492260 | 0.2234492160 | 0.2234495335 | 0.2234495335 |
| $x = 0.3$ | 0.62514 | 0.6251179353 | 0.6251179207 | 0.6251182341 | 0.6251182333 |
| $x = 0.5$ | 0.87728 | 0.8772800857 | 0.8772800822 | 0.8772796533 | 0.8772796530 |
| $x = 0.7$ | 0.83686 | 0.8369228762 | 0.8369228771 | 0.8369225592 | 0.8369225599 |
| $x = 0.9$ | 0.36573 | 0.36575425470 | 0.3657542612 | 0.3657544560 | 0.3657544557 |

| $\|E\|_\infty$ | 5.57E-07 | 8.97E-10 |

**Problem 2 [28]**

Consider the two-dimensional form of the Burgers equation (5) with the following exact solution [28]

$$u(x,y,t) = \frac{1}{1+e^{(x+y-t)/\varepsilon}} \tag{56}$$

where the computational domain is $D = \{(x,y,t) \in R^2 \times R | \ 0 \leq x,y \leq 1 \text{ and } t \in [0, t_f]\}$. The initial and boundary conditions can be taken from the exact solution (56).

The direction parameter $\theta$ plays a critical role in terms of both accuracy and stability. Exponential convergence of the current algorithm for both $\varepsilon = 1$ and $\varepsilon = 0.1$ is illustrated for various $\theta$ values as seen in Figure 3. We observe that all our theoretical error expectations are satisfied for all values of the direction parameter $\theta$. Natural solution behavior (57) drops rapidly from one to zero when the viscosity constant approaches zero. Thus, many numerical solvers fail to produce non-oscillatory results in the regions in which the rapid changes occur. In Figure 4, it has been shown that the present ChCM-IELDTM technique has ability to accurately capture the shock behavior for both $\varepsilon = 0.01$ and $\varepsilon = 0.005$. In Table 2, our motivation is to compare the current IELDTM integration method with the MATLAB solvers, *ode45* and *ode15s*, in terms of both accuracy and computational cost. As seen in the table, the present IELDTM-ChCM provides optimal results with less number of time steps than the rival schemes.

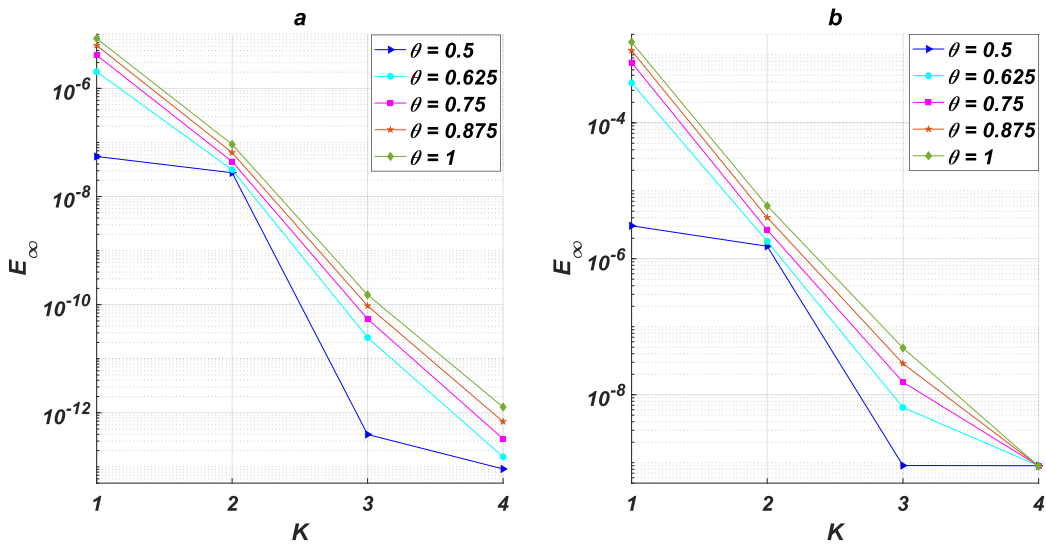

**Fig. 3** Order refinement results of the ChCM-IELDTM at $t_f = 0.5$ for various values of the direction parameter $\theta$ with the parameter values a) $\varepsilon = 1, N = M = 10$ and $dt = 0.005$ b) $\varepsilon = 0.1, N = M = 20$ and $dt = 0.005$.

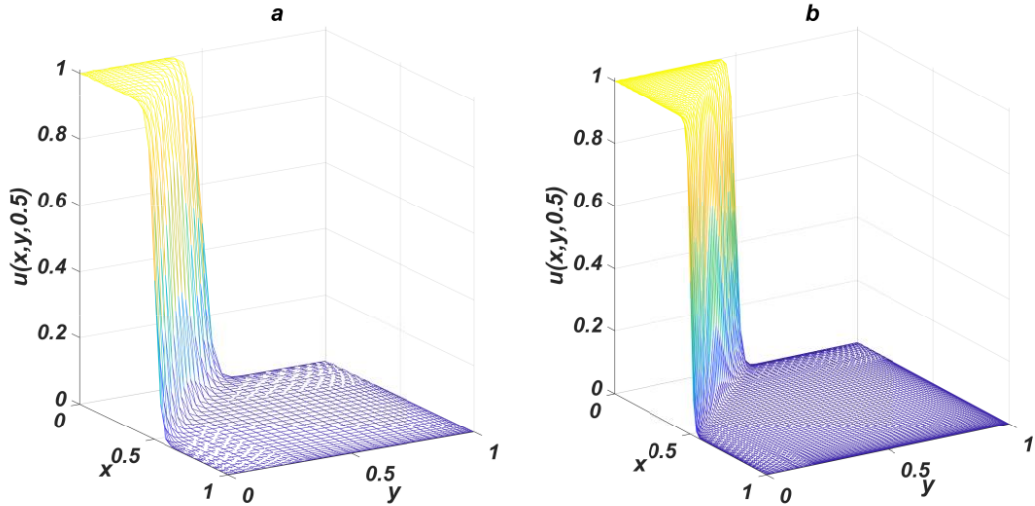

**Fig. 4** The ChCM-IELDTM solutions of Problem 2 at various times for the parameter values a) $\varepsilon = 0.01, N = M = 40, \theta = 0.5, K = 4$ and $dt = 0.01$ b) $\varepsilon = 0.005, N = M = 80, \theta = 0.5, K = 4$ and $dt = 0.01$.

**Table 2** Number of time steps and the maximum pointwise errors of the time-adaptive ChCM-IELDTM method for various values of the viscosity constant $\varepsilon$ and time approximation order $K$ with $\theta = 0.5$, $t_f = 0.5$ and $tol = E - 15$.

| $\varepsilon$ | Method | Step | $\|E\|_\infty$ |
|---|---|---|---|
| $\varepsilon = 1$ $N = M = 10$ | IELDTM-ChCSM $K = 3$ | 222 | 1.84E-14 |
| | IELDTM-ChCSM $K = 4$ | 178 | 1.48E-14 |
| | IELDTM-ChCSM $K = 5$ | 149 | 2.35E-14 |
| | *ode45*-ChCSM | 7023 | 1.87E-14 |
| | *ode15s*-ChCSM | 94 | 2.25E-14 |
| $\varepsilon = 0.1$ $N = M = 20$ | IELDTM-ChCSM $K = 3$ | 274 | 8.93E-10 |
| | IELDTM-ChCSM $K = 4$ | 220 | 9.02E-10 |
| | IELDTM-ChCSM $K = 5$ | 183 | 9.02E-10 |
| | *ode45*-ChCSM | 19919 | 9.02E-10 |
| | *ode15s*-ChCSM | 504 | 9.02E-10 |
| $\varepsilon = 0.01$ $N = M = 40$ | IELDTM-ChCSM $K = 3$ | 378 | 2.50E-03 |
| | IELDTM-ChCSM $K = 4$ | 311 | 2.60E-03 |
| | IELDTM-ChCSM $K = 5$ | 268 | 2.60E-03 |
| | ode45-ChCSM | 55436 | 2.60E-03 |
| | ode15s-ChCSM | 4327 | 2.60E-03 |

## 5. Conclusions and Recommendation

In this study, one-step implicit-explicit local differential transform method (IELDTM) has been derived for the temporal integration of problems represented by nonlinear advection-diffusion processes. The Chebyshev spectral collocation method (ChCM) has been considered for the spatial part of the considered PDEs due to computational efficiency of the method. The global error analysis of the IELDTM together with the stability analysis has been utilized and adaptivity equations have been proposed. The exponential convergence of the current numerical algorithm has been proven both theoretically and experimentally. The central case of the IELDTM, with the selection of $\theta = 0.5$, has been observed to be both more accurate and stability preserved. It has been proven that the IELDTM overcomes the shortcomings of current numerical techniques, such as the inaccuracy disadvantages of the θ-method and the instability drawbacks of the DTM-based methods. It has also been presented that the current time integration scheme achieves optimum accuracy by using fewer time steps than the commonly used time integrators of MATLAB such as *ode45* and *ode15s*. It has also been shown here that the ChCM-IELDTM has ability to solve advection dominated problems without producing any unwanted oscillation in both one and two dimensions.

### Acknowledgments

The first author would like to thank the Science Fellowships and Grant Programs Department of TUBITAK (TUBITAK BIDEB) for their support to his academic research.